\documentclass{elsart}

 \usepackage{graphicx}
 \usepackage{epsfig}
\usepackage{amssymb}

 \usepackage{lineno}

\usepackage{amsmath,amssymb,amsfonts,
latexsym, amscd, amsfonts, epsfig, eepic,epic}
\usepackage{mathrsfs}
\usepackage{color}
\usepackage{eucal}
\usepackage{eufrak}
\usepackage[all]{xypic}
\usepackage{xspace}
\theoremstyle{plain}

\newtheorem{theorem}{Theorem}[section]
\newtheorem{lemma}[theorem]{Lemma}

\theoremstyle{definition}

\theoremstyle{example}

\theoremstyle{remark}
\newcommand{\pff}{{\it Proof.}\ }

\begin{document}
\begin{frontmatter}



\title{Stacks in canonical RNA pseudoknot structures}


\author{Hillary S.W. Han and Christian M. Reidys$^{*}$}
\address{Center for Combinatorics, LPMC-TJKLC 
           \\
         Nankai University  \\
         Tianjin 300071\\
         P.R.~China\\
         Phone: *86-22-2350-6800\\
         Fax:   *86-22-2350-9272\\
         reidys@nankai.edu.cn}
\begin{abstract}
In this paper we study the distribution of stacks in $k$-noncrossing, 
$\tau$-canonical RNA pseudoknot structures ($\langle k,\tau\rangle
$-structures). 
An RNA structure is called $k$-noncrossing if it has no more than 
$k-1$ mutually crossing arcs and $\tau$-canonical if each
arc is contained in a stack of length at least $\tau$. Based on the
ordinary generating function of $\langle k,\tau\rangle$-structures 
\cite{Reidys:08ma} we derive the bivariate generating function 
${\bf T}_{k,\tau}(x,u)=\sum_{n \geq 0} \sum_{0\leq t \leq
\frac{n}{2}} {\sf T}_{k, \tau}^{} (n,t)\,u^t x^n$,
where ${\sf T}_{k,\tau}(n,t)$ is the number of 
$\langle k,\tau\rangle$-structures having exactly $t$ stacks
and study its singularities. We show that for a certain
parametrization of the variable $u$, ${\bf T}_{k,\tau}(x,u)$ has a
unique, dominant singularity. The particular shift of this
singularity parametrized by $u$ implies a central 
limit theorem for the distribution of stack-numbers.
Our results are of importance for understanding the ``language'' of 
minimum-free energy RNA
pseudoknot structures, generated by computer folding algorithms.
\end{abstract}

\begin{keyword}
$k$-noncrossing RNA structure, pseudoknot, generating function,
singularity analysis, central limit theorem
\end{keyword}
\end{frontmatter}

\section{Introduction}\label{S:Intro}
\begin{figure}[ht]
\centerline{%
\epsfig{file=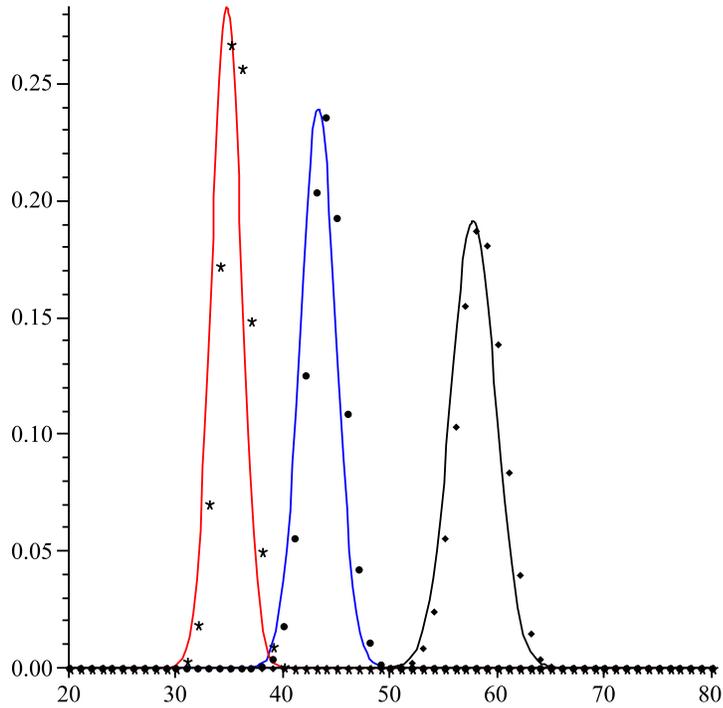,width=0.7\textwidth}\hskip15pt
 }
\caption{\small Central limit distributions for the stack-numbers in
$3$-noncrossing $\tau$-canonical RNA structures for $n=500$: we display 
the limit distribution of
Theorem~\ref{T:Gauss} for $\tau=3$ (black), $\tau=4$ (blue) and
$\tau=5$ (red). We also display the exact enumeration results obtained via
Theorem~\ref{T:exact}, eq.~(\ref{E:initial}) represented by $\diamond$,
$\bullet$ and $*$, respectively.}
\label{F:distribution}
\end{figure}

In this paper we compute the bivariate generating function
\begin{equation}
{\bf T}_{k,\tau}(x,u)= \sum_{n \geq 0} \sum_{0\leq t \leq
\frac{n}{2}} {\sf T}_{k, \tau}^{} (n,t)\,u^t x^n,
\end{equation}
where ${\sf T}_{k,\tau}(n,t)$ is the number of $k$-noncrossing,
$\tau$-canonical structures having exactly $t$ stacks.
Furthermore we prove the central limit theorem
\begin{equation}
\lim_{n\to\infty}\mathbb{P}\left(\frac{X_{k,\tau}^n- \mu_{k,\tau}
n}{\sqrt{n\, \sigma_{k,\tau}^2}} < x \right)  =
\frac{1}{\sqrt{2\pi}}\int_{-\infty}^{x}\,e^{-\frac{1}{2}t^2} dt \ ,
\end{equation}
where $X_{k,\tau}^n$ is a random variable with distribution
\begin{equation}\label{E:X}
\mathbb{P}(X_{k,\tau}^n=t)= {\sf T}_{k,\tau}(n,t)/{\sf
T}_{k,\tau}(n),
\end{equation}
see Fig.~\ref{F:distribution} and Tab.~\ref{T:1}. Our results allow to
compute the stack-numbers of a random $k$-noncrossing, $\tau$-canonical
RNA structure. Since the latter are obtained via the reflection principle,
it is nontrivial and currently not known how to construct such structures 
with uniform probability in polynomial time.

An RNA molecule is a sequence of the four nucleotides {\bf A}, {\bf
G}, {\bf U} and {\bf C} together with the Watson-Crick ({\bf A-U},
{\bf G-C}) and ({\bf U-G}) base pairing rules. RNA molecules form
``helical'' structures by pairing nucleotides and thereby
lowering their minimum free energy (mfe). The biochemistry of these
nucleotide-pairings favors parallel stacking of bonds due to
mesomerie effects. The resulting $3$-dimensional configuration of
the nucleotides is the RNA tertiary structure which determines the
functionality of the molecule.

\begin{figure}[ht]
\centerline{%
\epsfig{file=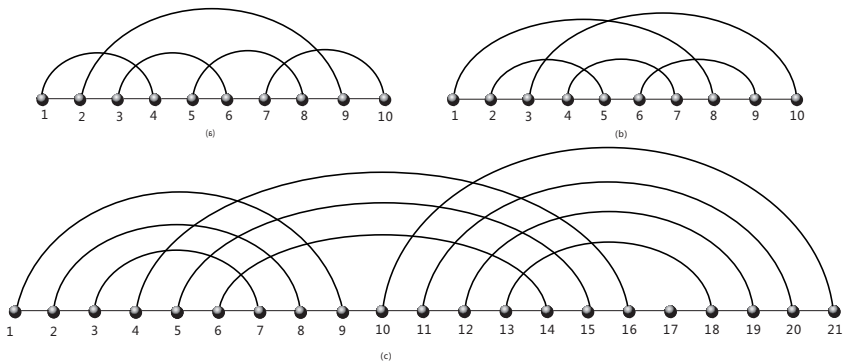,width=0.8\textwidth}\hskip15pt
 }
\caption{\small $k$-noncrossing diagrams: we display the only two nonplanar
$3$-noncrossing diagrams over $10$ vertices (upper) and a $3$-noncrossing,
$3$-canonical diagram with arc-length $\lambda\ge 4$ (lower). } \label{F:stack1}
\end{figure}

In this paper we study the distribution of stacks in
$k$-noncrossing, $\tau$-canonical RNA structures. We represent RNA
structures as diagrams and identify Watson-Crick ({\bf A-U}, {\bf
G-C}) and ({\bf U-G}) base pairings with arcs drawn in the upper
halfplane, ignoring the bonds of the primary sequence. A diagram is
a graph over the vertex set $[n]=\{1, \dots, n\}$ in which each
vertex has degree less or equal to one. It is represented by drawing
the vertices in a horizontal line and its arcs $(i,j)$, where $i<j$, in the
upper half-plane, see Fig.~\ref{F:stack1}.
The vertices and arcs correspond to nucleotides
and Watson-Crick ({\bf A-U}, {\bf G-C}) and ({\bf U-G}) base pairs,
respectively. Diagrams have the three key parameters $k$, $\lambda$
and $\tau$. Here $k-1$ is the maximum number of mutually crossing
arcs, $\lambda$ is the minimum arc-length and $\tau$ the minimum
length of a stack. A $\lambda$-arc is an arc $(i,j)$, where
$j-i=\lambda$ and a stack of length $\tau$ is a sequence of
``parallel'' arcs: $((i,j),(i+1,j-1),\dots,(i+(\tau-1),
j-(\tau-1)))$.

\begin{figure}[ht]
\centerline{%
\epsfig{file=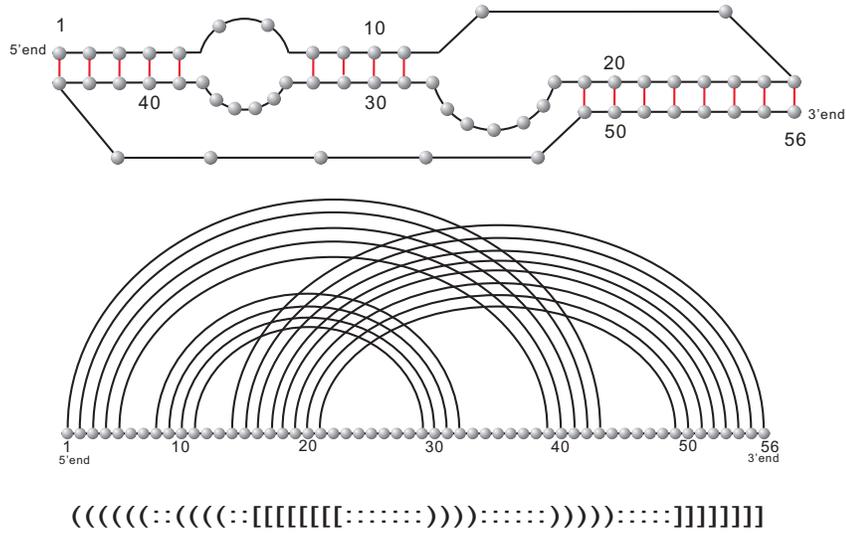,width=0.8\textwidth}\hskip15pt
 }
\caption{\small A $3$-noncrossing $3$-canonical RNA 
pseudoknot structures over $56$ nucleotides drawn as a planar graph (upper)
and as a $3$-noncrossing diagram (lower). } \label{F:stack2}
\end{figure}

We call a $k$-noncrossing diagram with arc-length $\lambda\ge 4$ and
stack-length $\tau\ge 3$ a $k$-noncrossing, $\tau$-canonical RNA
structure ($\langle k,\tau\rangle$-structure). 
We denote the number of $\langle k,\tau\rangle$-structures
and those with exactly $t$ stacks by ${\sf T}_{k,\tau}(n)$ and 
${\sf T}_{k,\tau} (n,t)$, respectively. For
$k=2$ this recovers the well known RNA secondary structures
\cite{Waterman:80,Waterman:78a,Waterman:79a,Waterman:94a}. Of
course, the interesting cases are $k\ge 3$, where we allow for
crossings, i.e.~the RNA pseudoknot structures. However, to the best
of our knowledge, our results are--even for RNA secondary structures
new. In Fig.~\ref{F:stack2} we give two representations of a typical
$\langle 3,3\langle$  RNA pseudoknot structure.

The results of this paper are relevant for the
understanding of mfe RNA pseudoknot structures. Computer folding
algorithms, like for instance, the dynamic programming routine of 
Rivas and Eddy \cite{Rivas:99a}, can generate, at least for 
$n\le 120$, mfe-pseudoknot structures
reasonably fast. 
While it is well known how to generate a random RNA secondary structure,
the particular construction hinges on the fact that the latter
can be build inductively. Although $\langle k,\tau\rangle$-structures are 
$D$-finite non inductive recursion exists \cite{Reidys:07pseu, Reidys:07lego} 
and it is at present time not known how to generate them with uniform 
probability. 
Our findings are therefore of particular importance for understanding the 
``language of RNA'', generated by {\it ab initio} folding algorithms. 
\section{Some basic facts}\label{S:basic}


In the following, we shall identify RNA pseudoknot structures with
$k$-non\-crossing diagrams with minimum arc-length $\ge 4$ and
stack-length $\tau\ge 3$. From now on we will always assume that any
structure has a minimum arc-length greater than four. 
Let ${T}^{}_{k,\tau}(n)$
(${\sf T}_{k,\tau}^{}(n)$) denote the set (number) of
$\langle k,\tau\rangle$-structures over $[n]$. 
Furthermore, let ${T}_{k,\tau}^{}(n,t)$ and ${\sf
T}_{k,\tau}^{}(n,t)$ denote the set (number) of $\langle k,
\tau\rangle$-structures having exactly $t$ stacks. In
addition, let ${T}_{k,\tau}^{}(n,h)$ and ${\sf T}_{k,\tau}^{}(n,h)$
denote the set (number) of $\langle k,\tau\rangle$-structures
having exactly $h$ arcs. The generating function,
${\bf T}_{k,\tau}^{}(x)=\sum_{n\ge 0}{\sf T}_{k,\tau}^{}(n)x^n$,
where $k\ge 2,\,\tau\ge 3$ of $\langle k,\tau\rangle$-structures has been
obtained in \cite{Reidys:08ma}. The function is closely related to
${\bf F}_k(x)=\sum_{n}f_k(2n)x^{2n}$, the generating function of
$k$-noncrossing matchings \cite{Chen}. Beyond functional equations implied
directly by the reflection-principle \cite{Grabiner:93a}, the
following asymptotic formula \cite{Wang:07} will be of importance
\begin{equation}\label{E:f-k-imp}
\forall\, k\in\mathbb{N};\qquad  f_{k}(2n) \, \sim  \, c_k  \,
n^{-((k-1)^2+(k-1)/2)}\, (2(k-1))^{2n},\qquad c_k>0 \ .
\end{equation}
${\bf T}^{}_{k,\tau}(x)$ is given
as follows:
\begin{theorem}\label{T:arc-4}
Let $k\ge 2$, $\tau\ge 3$, $x$ be an indeterminate and $\rho_k$ the
dominant, positive real singularity of ${\bf F}_k(z)$. Then ${\bf
T}^{}_{ k,\tau}(x)$ is given by
\begin{equation}\label{E:functional-arc-3}
{\bf T}_{k,\tau}^{}(x) = \frac{1}{v_0(x)} \ {\bf F}_k
\left(\frac{\sqrt{w_{0}(x)}\,x}{v_0(x)}\right) ,
\end{equation}
where $w_0(x) =\frac{x^{2\tau-2}}{1-x^2+x^{2\tau}}$ and $v_0(x)  =
1-x+w_0(x)x^2+w_0(x)x^3+w_0(x)x^4$. Furthermore,
\begin{equation}\label{E:growth-3}
{\sf T}_{k,\tau}^{}(n) \sim c_k\, n^{-(k-1)^2-\frac{k-1}{2}}\,
\left(\frac{1}{\gamma_{k,\tau}^{}}\right)^n\ ,\quad \mbox {\it
for}\quad k=2,3,\dots,9
\end{equation}
holds, where $\gamma_{k,\tau}^{}$ is the minimal positive real
solution of the equation $\frac{\sqrt{w_{0}(x)}\,x}{v_0(x)} =
\rho_k$.
\end{theorem}

In addition we shall make use of the following basic theorem on limit 
distributions:
\begin{theorem}\label{T:K}{\bf (L\'{e}vy-Cram\'{e}r)}
Let $\{\xi_{n}\}$ be a sequence of random variables and let
$\{\varphi_{n}(x)\}$ and $\{F_{n}(x)\}$ be the corresponding
sequences of characteristic and distribution functions. If there
exists a function $\varphi(t)$, such that $\lim_{n \rightarrow
\infty}\varphi_{n}(t)=\varphi(t)$ uniformly over an arbitrary finite
interval enclosing the origin, then there exists a random variable
$\xi$ with distribution function $F(x)$ such that
$$F_{n}(x)\Longrightarrow F(x)$$
uniformly over any finite or infinite interval of continuity of
$F(x)$.
\end{theorem}

\section{Combinatorics of stacks}\label{S:2}

\begin{figure}[ht]
\centerline{%
\epsfig{file=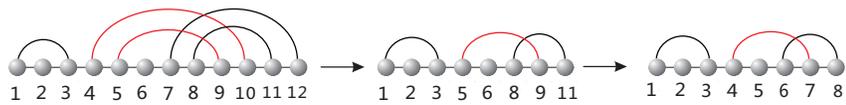,width=0.8\textwidth}\hskip15pt
 }
\caption{\small The mapping $c$: $T_{k,\tau}^{}(n,t)\rightarrow
\dot{\bigcup}_{t \tau \leq h \leq \lfloor \frac{n}{2}
\rfloor}\,C_k^{*}(n-2(h-t),t)$ is obtained in two steps: first
contraction of the stacks and second relabeling the resulting
diagram.} \label{F:core2}
\end{figure}

The objective of this section is to compute the bivariate
generating function
\begin{equation}
{\bf T}_{k,\tau}(x,u)= \sum_{n \geq 0} \sum_{0\leq t \leq
\frac{n}{2}} {\sf T}_{k, \tau}^{} (n,t)\,u^t x^n \ .
\end{equation}
For this purpose let us recall the concept of $k$-noncrossing cores
developed in \cite{Reidys:07lego}. A $k$-noncrossing core is a
$k$-noncrossing diagram in which all stacks have size one. We denote
the set and the number of cores over $[n]$ by $C_k(n)$ and ${\sf
C}_k(n)$, respectively. Furthermore, let $C_k(n,t)$  and ${\sf
C}_k(n,t)$ denote the set and number of cores having exactly $t$
stacks. We consider the arc-sets
$$
\beta_2=\{(i,i+2)\,| \ i+1 \ \mbox{\rm is isolated} \,\} \
\mbox{\rm and} \ \beta_3=\{(i,i+3)\,|\ i+1,i+2 \ \mbox{\rm are
isolated} \,\}$$ and set $\beta= \beta_2 \cup \beta_3 $.
Furthermore, let
\begin{eqnarray}
C_k^{*}(n,t)& = & \{\zeta \,|\, \zeta \in C_k(n,t);\,
\zeta \ \mbox{\rm contains no $\beta$-arc} \, \} \\
T_{k,1}^{*}(n,h) & = & \{\zeta \,|\, \zeta \in T_{k,1}(n,h);\, \zeta
\ \mbox{\rm contains no $\beta$-arc} \, \}
\end{eqnarray}
where $T_{k,1}(n,h)$ denotes the set of $k$-noncrossing diagrams without
$1$-arcs having exactly $h$ arcs. 
We set ${\sf C}_k^{*}(n,0)={\sf T}_{k,1}^{*}(n,0)=1$ for $0\le
n$. 
In our first theorem we shall show that the number of all
$\langle k,\tau\rangle$-structures with exactly $t$ stacks is a sum
of the number of $C_k^*(n',t)$-diagrams with positive integer
coefficients. Let $f_k(m,\ell)$ denote the number of $k$-noncrossing 
partial matchings having exactly $\ell$ isolated points. 
\begin{theorem}\label{T:exact}
Suppose we have $k,t, \tau \in \mathbb{N}$, $k\ge 2$, $\tau \geq 3$
and $1\leq t \leq \lfloor n/2\rfloor$. Then
\begin{equation}\label{E:initial}
{\sf T}_{k,\tau}^{}(n,t) = \sum_{h \geq t \tau}^{\lfloor n/2\rfloor}
{{(1-\tau)t+h-1}\choose{t-1}}\,{\sf C}_k^{*}(n+2t-2h,t) \\
\end{equation}
\begin{equation}\label{E:inversion}
\begin{split}
{\sf C}_k^{*}(n,h) = & \sum_{b=0}^{h-1}
(-1)^{h-b-1} {h-1 \choose b}\, {\sf T}_{k,1}^{*}(n-2h+2b+2,b+1)\\
& \quad for \quad h \geq 1.
\end{split}
\end{equation}
Furthermore, ${\sf T}_{k,1}^{*}(n,h)$ satisfies
\begin{equation}\label{E:kl}
\begin{split}
{\sf T}_{k,1}^{*}(n,h)=&\sum_{0 \leq j_1+j_2+j_3 \leq h} (-1)^{j_1+j_2+j_3}\\
& \lambda(n, j_1,j_2,j_3)\, f_k(n-2j_1-3j_2-4j_3, n-2h-j_2-2j_3)
\end{split}
\end{equation}
where
$$
\lambda(n, j_1,j_2,j_3)= {n-j_1-2j_2-3j_3 \choose j_1,j_2,j_3,
n-2j_1-3j_2-4j_3} .
$$
\end{theorem}
\pff First, there exists a mapping from $\langle k,\tau\rangle$-structures
with $t$ stacks into
$C_k^{*}(n-2(h-t),t)$:
\begin{equation}
c\colon T_{k,\tau}^{}(n,t)\rightarrow \dot{\bigcup_{t \tau \leq h
\leq \lfloor \frac{n}{2} \rfloor}}\,C_k^{*}(n-2(h-t),t), \qquad\
\zeta \mapsto c(\zeta)
\end{equation}
The core, $c(\zeta)$, is obtained in two steps:
first we map arcs and isolated vertices as follows
\begin{equation}
\forall \ell \geq \tau-1; \quad ((i-\ell, j+\ell),\ldots, (i,j))
\mapsto (i,j) \ \text{\rm and} \ j \mapsto j \quad
\text{\rm if $j$ is isolated.}
\end{equation}
Second we relabel the vertices of the resulting diagram from left to
right in increasing order. That is we replace each stack by a single
arc, keep isolated vertices and then relabel, see Fig.~\ref{F:core2}.
It is straightforward to verify that $c\colon T_{k,\tau}^{}(n,t)\rightarrow
\dot{\bigcup}_{t \tau \leq h\leq \lfloor \frac{n}{2}
\rfloor}\,C_k^{*}(n-2(h-t),t)$ is well-defined and surjective.
Taking into account that each stack has its specific length gives
rise to consider
\begin{equation}\label{E:f}
\begin{split}
f_{k,\tau}\colon   T_{k,\tau}^{}(n,t)  &\longrightarrow
\dot{\bigcup_{ t\,\tau \leq h \leq \lfloor
\frac{n}{2}\rfloor}} (C_k^{*}(n+2t-2h,t)\times  \\
&  \{(\alpha_i)_{1 \leq i \leq
t}\mid\sum_{i=1}^{t}\alpha_i=h-t,\alpha_i \geq \tau-1\})
\end{split}
\end{equation}
given by $f_{k,\tau}(\zeta)=(c(\zeta),(\alpha_i)_{1 \leq i \leq t
})$. By construction $f_{k,\tau}$ is well-defined and a bijection:
obviously we can reconstruct the original
$T_{k,\tau}^{}(n,t)$-element from its core and the sequence of its
stack-multiplicities (labeling the stacks from left to right).
Computing the multiplicities of the resulting cores we derive
\begin{equation}
\left| \left\{(\alpha_i)_{1 \leq i \leq t}\mid \sum_{i=1}^{t} \alpha_i=h-t;
\alpha_i \geq \tau-1\right\}\right|={(1-\tau)t+h-1 \choose t-1} \ .
\end{equation}
We can therefore conclude
\begin{equation}
{\sf T}_{k,\tau}^{}(n,t)= \sum_{h \geq t \tau}^{\lfloor n/2 \rfloor}
 {{(1-\tau)t+h-1}\choose{t-1}}\,{\sf C_k}^{*}(n+2t-2h,t),
\end{equation}
whence eq.~(\ref{E:initial}). Next we switch and consider arcs
instead of stacks, using the parameter $h$. Again, collapsing the
stack induces the mapping
\begin{equation}
c^{*}\colon  T_{k,1}^{*}(n,h)\rightarrow \dot{\bigcup}_{0 \leq b
\leq h-1 }\,C_k^{*}(n-2b,h-b)  \qquad\ \zeta \mapsto c^{*}(\zeta)
\end{equation}
Indeed, $c^{*}$ is well defined, since a $T_{k,1}^{*}(n,h)$-diagram
can be mapped into a core structure without $1$- and $\beta$- arcs,
i.e.~into an element of $C_k^{*}(n',h')$. It is straightforward to
show that $c^*$ is surjective. Labeling the $h-b$ stacks of $\zeta
\in T_{k,1}^{*}(n,h)$ from left to right and keeping track of
multiplicities gives rise to the map
\begin{equation}\label{E:m}
\begin{split}
m_{k}\colon  T_{k,1}^{*}(n,h) & \rightarrow \\
&
\dot{\bigcup_{0 \leq b\leq h-1}}\left[ C_k^{*}(n-2b,h-b)\times
\left\{(a_i)_{1 \leq i \leq h-b}\, |\, \sum_{i=1}^{h-b}a_i=b, a_i \geq
0\right\}\right]
\end{split}
\end{equation}
given by $m_{k}(\zeta)=(c^*(\zeta),(a_i)_{1 \leq i \leq h-b})$. The
mapping $m_{k}$ is well-defined and a bijection. Clearly,
\begin{equation}\label{E:bino}
\left|\left\{(a_i)_{1 \leq i \leq h-b}\mid \sum_{i=1}^{h-b} a_i=b; a_i \geq
0 \right\} \right|={ h-1 \choose b}
\end{equation}
and eq.~(\ref{E:m}) and eq.~(\ref{E:bino}) imply
\begin{equation}
{\sf T}_{k,1}^{*}(n,h) = \sum_{b=0}^{h-1}{ h-1 \choose b}\,{\sf
C}_k^{*}(n-2b,h-b).
\end{equation}
Via M\"obius-inversion eq.~(\ref{E:inversion}) follows. It is
straightforward to show that there are
$$
\lambda(n, j_1,j_2,j_3) = {n-j_1-2j_2-3j_3 \choose j_1, j_2, j_3,
n-2j_1-3j_2-4j_3}
$$
ways to select $j_1$ $1$-arcs, $j_2$ $\beta_2$-arcs and $j_3$
$\beta_3$-arcs over $[n]$. Since removing $j_1$ $1$-arcs, $j_2$
$\beta_2$-arcs and $j_3$ $\beta_3$-arcs removes $2j_1+3j_2+4j_3$
vertices, the number of configurations of at least $j_1$ $1$-arcs,
$j_2$ $\beta_2$-arcs and $j_3$ $\beta_3$-arcs is given by
$\lambda(n, j_1, j_2, j_3) f_k(n-2j_1-3j_2-4j_3, n-2h-j_2-2j_3)$.
Via the inclusion-exclusion principle, we arrive at
\begin{equation}
\begin{split}
{\sf T}_{k,1}^{*}(n,h) &= \sum_{ 0 \leq j_1+j_2+j_3 \leq h}(-1)^{j_1+j_2+j_3}\\
& \lambda (n,j_1,j_2,j_3) \, f_k(n-2j_1-3j_2-4j_3, n-2h-j_2-2j_3),
\end{split}
\end{equation}
whence Theorem~\ref{T:exact}.\qed

Of course, Theorem~\ref{T:exact} implies a functional relation for
${\bf T}_{k,\tau}(x,u)$:
\begin{lemma}\label{L:func}
Let $k\ge 2$, $\tau \ge 3$ and let u, x be indeterminates. Then we
have the functional relation
\begin{equation}
\sum_{n \geq 0} \sum _{0 \leq t \leq \lfloor\frac{n}{2}\rfloor} {\sf
T}_{k,\tau}^{}(n,t)u^t x^n  = \sum_{n \geq 0} \sum _{0 \leq t \leq
\lfloor\frac{n}{2}\rfloor} {\sf C}_k^*(n,t)\left(\frac{u
x^{2(\tau-1)}}{1-x^2}\right)^t x^n
\end{equation}
and in particular, for $u=1$
\begin{equation}
\sum_{n \geq 0} {\sf T}_{k,\tau}^{}(n) x^n  = \sum_{n \geq 0} \sum
_{0\leq t \leq \lfloor\frac{n}{2}\rfloor} {\sf
C}_k^*(n,t)\left(\frac{ x^{2(\tau-1)}}{1-x^2}\right)^t x^n.
\end{equation}
\end{lemma}

\pff We set
$$
\sum _{t\ge 1} \left[ \sum_{n \geq 2t\tau} 
{\sf C}_k^*(n,t) x^n\right]u^t = \sum_{t \geq 1} \varphi_t(x)u^t
$$
and proceed by deducing the functional equation for 
${\bf T}_{k,\tau}(x,u)$ via Theorem~\ref{T:exact}. For this purpose we
note that for $t=0$ the Binomial coefficient
$
{(1-\tau)t+h-1\choose t-1} 
$
is zero, while the term ${\sf T}_{k,\tau}^{}(n,0)=1$ for $n \geq 1$.
Clearly, ${\sf T}_{k,\tau}^{}(n,0)=1$ counts for each $n \geq 1$ the
structure consisting only of isolated vertices. We accordingly
have to extend the identity of Theorem~\ref{T:exact}
$$
{\sf T}_{k,\tau}^{}(n,t)= \sum_{h \geq t \tau}^{\lfloor n/2 \rfloor}
 {{(1-\tau)t+h-1}\choose{t-1}}\,{\sf C}_k^{*}(n+2t-2h,t)
$$
to the case $t=0$, $n \geq 1$ which gives rise to the term $\sum_{n
\geq 1}x^n=\frac{x}{1-x}$. Accordingly, we derive
\begin{equation}\label{E:generating}
\begin{split}
\sum_{n \geq 0} \sum _{0\leq t \leq \lfloor\frac{n}{2}\rfloor}{\sf
T}_{k,\tau}^{}(n,t)u^t x^n= &  
\sum_{n \geq 2} \sum _{1 \leq t \leq
\lfloor\frac{n}{2}\rfloor}\sum_{h \geq t \tau}^{\lfloor n/2
\rfloor}
{\sf C}_k^{*}(n+2t-2h,t) \; \times \\
& \qquad \qquad  {{(1-\tau)t+h-1}\choose{t-1}}u^t x^n
+\sum_{n \geq 0}x^n.
\end{split}
\end{equation}
We rewrite the right hand side of eq.~(\ref{E:generating})
\begin{align*}
=& \sum_{t \geq 1}\sum_{h \geq t \tau}^{\lfloor n/2 \rfloor}\sum_{n
\geq 2t\tau}{\sf
C}_k^{*}(n+2t-2h,t)x^{n+2t-2h}{{(1-\tau)t+h-1}\choose{t-1}}u^t
x^{2h-2t} + \frac{1}{1-x}\\
=&\sum_{t \geq 1}\sum_{h \geq t \tau}^{\lfloor n/2
\rfloor}\varphi_t(x){{(1-\tau)t+h-1}\choose{t-1}}u^t x^{2h-2t} +
\frac{1}{1-x}.
\end{align*}
Rearranging the terms of the summation we obtain
\begin{align*}
=&\sum_{t \geq 1}\sum_{h \geq t \tau}^{\lfloor n/2
\rfloor}\varphi_t(x){{t-t \tau +h-1}\choose{h-t \tau}}x^{2h}
(\frac{u}{x^2})^t + \frac{1}{1-x}\\
=&\sum_{t \geq 1}\sum_{h \geq t \tau}^{\lfloor n/2 \rfloor}
\varphi_t(x){{h-t \tau +(t-1)}\choose{h-t \tau}}(x^2)^{h-t \tau}
(u x^{2(\tau-1)})^t + \frac{1}{1-x} .
\end{align*}
Using $\sum_n {{n+k} \choose n}x^n =\frac{1}{(1-x)^{k+1}}$, we can
transform the summation over $h$ and derive
\begin{align*}
\sum_{n \geq 0} \sum _{0\leq t \leq \lfloor\frac{n}{2}\rfloor}{\sf
T}_{k,\tau}^{}(n,t)u^t x^n =& \sum_{t \geq 1}\varphi_t(x)
\left(\frac{1}{1-x^2}\right)^t (u
x^{2(\tau-1)})^t+ \frac{1}{1-x}\\
=&\sum_{n \geq 0} \sum _{0 \leq t \leq \lfloor\frac{n}{2}\rfloor}
{\sf C}_k^*(n,t)\left(\frac{u x^{2(\tau-1)}}{1-x^2}\right)^t x^n  
\end{align*}
and the proof of Lemma~\ref{L:func} is complete.
\qed

We next consider a functional equation for $\sum_{n,h}{\sf
T}_k^{*}(n,h)u^hx^n$ proved in \cite{Reidys:07lego}.


\begin{lemma}\label{L:h-generate}
Let $k,\tau\in \mathbb{N}$, $k\ge 2$ and let $u,x$ be
indeterminates. Suppose we have
\begin{equation}\label{E:uni1}
\begin{split}
\forall\, h\ge 1,\ \ & {\sf A}_{k,\tau}^{}(n,h)  =\\
&\sum_{b=\tau-1}^{h-1} {b+(2-\tau)(h-b)-1 \choose h-b-1} {\sf
B}_k(n-2b,h-b)\ \text{\it and} \  {\sf A}_{k,\tau}^{}(n,0)=1 \ .
\end{split}
\end{equation}
Then we have the functional relation
\begin{eqnarray}\label{E:universal}
\sum_{n \ge 0}\sum_{0\le h\le \frac{n}{2}}{\sf A}^{ }_{k,
\tau}(n,h)u^hx^n & = & \sum_{n \ge 0}\sum_{0\le h \le
\frac{n}{2}}{\sf B}^{ }_{k} (n,h)\left(\frac{u\cdot (ux^2)^{
\tau-1}}{1-ux^2}\right)^hx^n\ .
\end{eqnarray}
\end{lemma}

Combining Lemma \ref{L:h-generate} and $ {\sf T}_{k,1}^{*}(n,h) =
\sum_{b=0}^{h-1}{ h-1 \choose b}\,{\sf C}_k^{*}(n-2b,h-b) $ we
arrive at
\begin{eqnarray}
\label{E:corestar} \sum_{n \ge 0}\sum_{0\le h\le \frac{n}{2}}{\sf
T}^{*}_{k,1}(n,h)u^h x^n  =  \sum_{n \ge 0}\sum_{0 \le h \le
\frac{n}{2}}{\sf C}^{
*}_{k}(n,h)\left(\frac{u}{1-ux^2}\right)^hx^n \ .
\end{eqnarray}
We shall make use of an additional relation for 
$\sum_{n \ge 0}\sum_{0\le h\le \frac{n}{2}}{\sf
T}^{*}_{k,1}(n,h)u^h x^n$, proved in \cite{Reidys:08ma}:
\begin{lemma}\label{L:func2}
Let $k \in \mathbb{N}$, $k \geq 2$, and x, w be indeterminates. Then
we have the functional relation
\begin{equation}\label{E:dagger}
\begin{split}
 \sum_{n\ge 0}\sum_{h \le \frac{n}{2}}{\sf T}_{k,1}^{*}(n,h)w^h x^n
=\frac{1}{v(x)}\sum_{n \ge
0}f_k(2n)\left(\frac{\sqrt{w}x}{v(x)}\right)^{2n} \ ,
\end{split}
\end{equation}
where $v(x) = 1-x+w x^2+w x^3+w x^4$ and $f_k(2n)$ is the number of 
$k$-noncrossing matchings over $2n$ vertices.
\end{lemma}
We are now in the position to present the main result of this section:
\begin{theorem}\label{T:f-eq}
Let $k\ge 2$, $\tau \ge 3$ and suppose $u$, $x$ are indeterminates.
Then we have the identity of formal power series
\begin{equation}\label{E:last}
\sum_{n \geq 0} \sum_{0\leq t \leq \frac{n}{2}} {\sf T}_{k, \tau}^{}
(n,t)u^t x^n= \frac{1}{v(x)}\sum_{n \ge
0}f_k(2n)\left(\frac{\sqrt{u_0}\,x}{v(x)}\right)^{2n},
\end{equation}
where $v(x)$ and $u_0= u_0(x,u)$ are given by
\begin{eqnarray}
u_0 & = & \frac{u x^{2(\tau-1)}}{u x^{2\tau}-x^2 +1}\\
v(x) & = & 1-x+u_0 x^2+u_0 x^3+u_0 x^4.
\end{eqnarray}
In particular, we can consider eq.~(\ref{E:last}) as a relation
between analytic functions, valid for $u=e^s$ and $|s| < \epsilon$
for $\epsilon$ sufficiently small and $|x|\le 1/2$.
\end{theorem}
\pff We interpret the bivariate generating function 
$$
\sum_{n \ge 0} \sum_{0 \le h \le \frac{n}{2}}{\sf C}^{*}_{k}(n,h)y^h x^n
$$ 
in two different ways: first, via Lemma~\ref{L:func} we have
\begin{equation}
\sum_{n \geq 0} \sum _{0 \leq t \leq \lfloor\frac{n}{2}\rfloor} {\sf
T}_{k,\tau}^{}(n,t)u^t x^n  = \sum_{n \geq 0} \sum _{0 \leq t \leq
\lfloor\frac{n}{2}\rfloor} {\sf C}_k^*(n,t)\left(\frac{u
x^{2(\tau-1)}}{1-x^2}\right)^t x^n
\end{equation}
and second according to eq.~(\ref{E:corestar}):
\begin{equation}\label{E:DD}
\sum_{n \ge 0}\sum_{0\le h\le \frac{n}{2}}{\sf
T}^{*}_{k,1}(n,h)u_0^h x^n  =  \sum_{n \ge 0} \sum_{0 \le h \le
\frac{n}{2}}{\sf C}^{
*}_{k}(n,h)\left(\frac{u_0}{1-u_0 x^2}\right)^h x^n.
\end{equation}
Note that, by definition of cores, we may replace the index ``$h$'' by ``$t$''
in the right hand side of eq.~(\ref{E:DD}).
The key observation is the relation between the terms $\frac{u
x^{2(\tau-1)}}{1-x^2}$ and $\frac{u_0}{1-u_0 x^2}$. Using the Ansatz
\begin{equation}
\frac{u x^{2(\tau-1)}}{1-x^2}= \frac{u_0}{1-u_0 x^2},
\end{equation}
we obtain the unique solution $u_0=\frac{ux^{2(\tau-1)}}{1-x^2+u
x^{2 \tau}}$. Accordingly, we conclude
\begin{align*}
\sum_{n \geq 0} \sum _{0 \leq t \leq \lfloor\frac{n}{2}\rfloor} {\sf
T}_{k,\tau}^{}(n,t)u^t x^n = & \sum_{n \geq 0} \sum _{0 \leq t \leq
\lfloor\frac{n}{2}\rfloor} {\sf C}_k^*(n,t)\left(\frac{u
x^{2(\tau-1)}}{1-x^2}\right)^t x^n \\
= & \sum_{n \ge 0}\sum_{0\le h\le \frac{n}{2}}{\sf
T}^{*}_{k,1}(n,h)u_0^h
x^n  \\
= & \frac{1}{v(x)}\sum_{n \ge
0}f_k(2n)\left(\frac{\sqrt{u_0}\,x}{v(x)}\right)^{2n},
\end{align*}
where $v(x) = 1-x+u_0 x^2+u_0 x^3+u_0 x^4$, and $u_0 =\frac{u
x^{2(\tau-1)}}{u x^{2\tau}-x^2 +1}$. In order to consider 
eq.~(\ref{E:last}) as a relation between analytic functions 
we need to satisfy $u
x^{2\tau}-x^2 +1 \neq 0$. Suppose $u=e^s$, $|s|<\epsilon$ for
$\epsilon$ sufficiently small. From $|x| \le  1/2$ and the
continuity (in $s$) of the roots of the family of polynomials
\begin{equation*}
\omega_s(X)=e^s X^{2\tau}-X^2 +1,\qquad \vert s\vert <\epsilon
\end{equation*}
we conclude $\omega_s(x) \neq 0$. Therefore eq.~(\ref{E:last}) holds
for $u=e^s$, $|s|<\epsilon$ and sufficiently small $\epsilon$ and
$|x| \le  1/2$.\qed

\section{The central limit theorem}
Suppose $\epsilon >0$, $2\le k$ and $u$ is parametrized as $u=e^s$,
where $|s|< \epsilon$. We set
\begin{eqnarray}
\varphi_{n,k,\tau}(s) & = &  \sum_{t \leq \frac{n}{2}} {\sf T}_{k,
\tau}(n, t) e^{ts} \\
{\bf U}_k(z,s) & = & \sum_{n \geq 0} \varphi_{n,k,\tau}(s) z^n .
\end{eqnarray}
We will use Theorem~\ref{T:f-eq}, which relates the generating functions 
${\bf U}_{k}(z,s)$ and ${\bf F}_k(z)=\sum_n f(2n)z^{2n}$ in order 
to study the singularities of ${\bf U}_k(z,s)$.
\begin{theorem}\label{T:z}
Suppose $\epsilon > 0$, $2\le k\le 9$, $3\le \tau\le 7$ and $u=e^s$,
where $|s|< \epsilon$. Then the following assertions hold:\\
{\sf (a)} For sufficiently small $\epsilon$, ${\bf U}_k(z,s)$ has
the unique singularity, $\gamma_{k,\tau}^{}(s)$, which is analytic
in $s$ and the unique minimal real positive solution of
\begin{equation}\label{E:singularity}
\frac{\sqrt{u_0(s)}\,z}{1-z+u_0(s) z^2+u_0(s) z^3+u_0(s)
z^4}-\rho_k=0.
\end{equation}
{\sf (b)} The coefficients of ${\bf U}_k(z,s)$ are asymptotically
given by
\begin{equation}\label{E:uniform}
[z^n] {\bf U}_k(z,s)= A(s)\, (1-O(n^{-1}))\
n^{-((k-1)^2+\frac{k-1}{2})} \left( \frac{1}{\gamma_{k,\tau}^{}(s)}
\right)^n, \quad \text{\it $A(s) \in \mathbb{C}$},
\end{equation}
uniformly in $s$ in a neighborhood of $0$. In particular, the
subexponential factors of the coefficients of ${\bf U}_k(z,s)$
coincide with those of ${\bf F}_k(z)$ and are independent of $s$.
\end{theorem}
\pff According to Theorem~\ref{T:f-eq} we have for $|s|< \epsilon$
\begin{equation}\label{E:u(s)}
{\bf U}_k(z,s)=\frac{1}{v(z,s)}\sum_{n \ge
0}f_k(2n)\left(\frac{\sqrt{u_0(s)}\,z}{v(z,s)}\right)^{2n} \
\end{equation}
where $v(z,s) = 1-z+u_0(s) z^2+u_0(s) z^3+u_0(s) z^4$ and $u_0(s)
=\frac{e^s z^{2(\tau-1)}}{e^s z^{2\tau}-z^2 +1}$. We set
\begin{eqnarray}
\psi_{\tau}(z,s) & = &
\frac{\sqrt{u_0(s)}\,z}{1-z+u_0(s) z^2+u_0(s) z^3+u_0(s) z^4} \\
{\bf W}_k(z,s) & = & \sum_{n \geq 0}f_k(2n)\left(\frac{
\sqrt{u_0(s)}\,z}{v(z,s)}\right)^{2n}.
\end{eqnarray}
Our first objective is to prove the existence of the unique
singularity, $\gamma_{k,\tau}^{}(s)$. For this purpose we consider
\begin{equation}\label{E:RR}
F(z,s)=\psi_{\tau}(z,s)- \rho_k .
\end{equation}
For $s=0$ and $2\le k\le 9$, $3\le \tau\le 7$ it is straightforward
to verify that a unique minimal real solution, $\omega$, exists.
We observe that for $|s|<\epsilon$ the following holds: 
(i) $F(\omega,0)=0$, (ii)
$F_z(\omega,0)\neq 0$ and (iii) the partial derivatives $F_z(z,s)$
and $F_s(z,s)$ are continuous. According to the analytic implicit
function theorem \cite{Flajolet:07a}, there exists in a sufficiently
small neighborhood of $0$, a unique analytic function
$\gamma_{k,\tau}^{}(s)$, that satisfies
\begin{equation*}
\forall\, s;\,\vert s\vert <\epsilon;\quad F(\gamma_{k,\tau}(s),s)=0
\quad \text{\rm and} \quad \gamma_{k,\tau}(0)=\omega.
\end{equation*}
{\it Claim $1$.} For $\epsilon$ sufficiently small and
$|s|<\epsilon$,
$\gamma_{k,\tau}^{}(s)$ is a dominant singularity of ${\bf U}_k(z,s)$.\\
Let $\zeta(s)$ be a dominant singularity of ${\bf U}_k(z,s)$. It is
straightforward to prove that $\zeta(0)$ is necessarily a singularity
of ${\bf W}_k(z,0)$. We proceed by applying
an continuity argument. For $\epsilon$ sufficiently small and $
\vert s\vert<\epsilon $ the singularities of $v(z,s)^{-1}$ and
$\gamma_{k,\tau}^{}(s)$ are both continuous in $s$. Therefore we can
conclude from our observation for $s=0$ that, for sufficiently small
$\epsilon$, all singularities of $v(z,s)^{-1}$ have modulus strictly
larger than $\gamma_{k,\tau}^{}(s)$. Accordingly, for sufficiently
small $\epsilon$, $\gamma_{k,\tau}^{}(s)$ is a dominant singularity
of
${\bf U}_k(z,s)$ , whence Claim $1$.\\
{\it Claim $2$.} $\gamma_{k,\tau}^{}(s)$ is unique.\\
Functional relations arising from the reflection principle
\cite{Grabiner:93a} imply that the generating function ${\bf
F}_k(z)$ is $D$-finite \cite{Stanley:80}. Accordingly there exists
some $e\in \mathbb{N}$ for which ${\bf F}_k(z)$ satisfies an ODE of
the form
\begin{equation}\label{E:JK}
q_{0,k}(z)\frac {d^e}{dz^e}{\bf
F}_k(z)+q_{1,k}(z)\frac{d^{e-1}}{dz^{e-1}}{\bf F}_k(z)+
q_{e,k}(z){\bf F}_k(z)=0 ,
\end{equation}
where $q_{j,k}(z)$ are polynomials. Any dominant singularity of
${\bf F}_k(z)$ is contained in the set of roots of $q_{0,k}(z)$
\cite{Stanley:80} and via direct computation we can verify that
$\gamma_{k,\tau}^{}(0)$ is the {\it unique} solution with minimal
modulus of all solutions of
\begin{equation}
\psi_{\tau}(z,0)=\vert \rho_k\vert .
\end{equation}
The analytic implicit function theorem applied to eq.~(\ref{E:RR})
now implies locally the existence of an unique analytic function
$\gamma_{k,\tau}^{}(s)$ solving $\psi_{\tau}(z,s)=\rho_k$. Using
continuity we conclude from the fact that $\gamma_{k,\tau}^{}(0)$ is
the {\it unique} solution with minimal modulus of all solutions of
$\psi_{\tau}(z,0)=\vert \rho_k\vert$, that for $\epsilon$
sufficiently small, $\gamma_{k,\tau}^{}(s)$ is the { unique}
solution with minimal modulus of all solutions of $\psi_{\tau}(z,s)=\vert
\rho_k\vert$.
This implies Claim $2$ and {\sf (a)} follows.\\
We finally prove {\sf (b)}. According to eq.~(\ref{E:f-k-imp}) we
have $$ f_{k}(2n) \, \sim  \, c_k  \, n^{-((k-1)^2+(k-1)/2)}\,
(2(k-1))^{2n}$$ for some $c_k>0$ and 
\begin{eqnarray*}
{\bf F}_k(z)=
\begin{cases}
O((z-\rho_k)^{(k-1)^2+(k-1)/2-1} \ln^{}(z-\rho_k)) & \text{\rm for
$k$ odd,
$z\rightarrow \rho_k$} \\
O((z-\rho_k)^{(k-1)^2+(k-1)/2-1})             & \text{\rm for $k$
even, $z\rightarrow \rho_k$,}
\end{cases}
\end{eqnarray*}
in accordance with basic structure theorems for singularities of
solutions of eq.~(\ref{E:JK}) \cite{Flajolet:07a}, p.~$499$. Let
$Q_{\gamma_{k,\tau}(s)}(z,s)$ denote the singular expansion of ${\bf
U}_k(z,s)$ at $\gamma_{k,\tau}(s)$. We have shown in Claim $1$ that
$\psi_{\tau}(z,s)$ does not
induce any dominant singularities and is regular at $\rho_k$. Let
$Q_{\rho_k}(z)$ denote the singular expansion of ${\bf F}_k(z)$ at
the dominant singularity $\rho_k$, i.e. $F_k(z) = O(Q_{\rho_k}(z))$
for $z\rightarrow \rho_k$. The singular expansion of the compositum,
${\bf F}_k(\psi_{\tau}(z,s))$, i.e.~$Q_{\gamma_{k,\tau}(s)}(z,s)$,
is derived by substituting the Taylor-expansion of
$\psi_{\tau}(z,s)$ into $Q_{\rho_k}(z)$ and we observe
\begin{equation}\label{E:si2}
Q_{\gamma_{k,\tau}(s)}(z,s)= Q_{\rho_k}(\psi_{\tau}(\zeta_k(s),s))
=O(Q_{\gamma_{k,\tau}(s)}(z))\ .
\end{equation}
Indeed, eq.~(\ref{E:si2}) follows immediately substituting
$\psi_{\tau}(z,s)-\psi_{\tau}(\gamma_{k,\tau}(s),s)$ for $z-\rho_k$
which does not change the singular expansion. Since
$\psi_{\tau}(z,s)$ is continuous in $s$ the singular expansion is
uniform with respect to the parameter $s$ and the uniformity lemma
of singularity analysis \cite{Flajolet:07a} implies
\begin{equation}\label{E:si3}
[z^n]\, {\bf U}_k(z,s) =  A(s)\,\left(1-O({1}/{n})\right)\,
n^{-((k-1)^2+(k-1)/2)}\,
\left(\frac{1}{\gamma_{k,\tau}(s)}\right)^n\quad 
\end{equation}
for some $A(s)\in\mathbb{C}$, uniformly in $s$ in a neighborhood of $0$. 
Therefore the asymptotic expansion is uniform in $s$ and 
eq.~(\ref{E:uniform}) follows. 
The proof shows in addition that the subexponential factors of the
coefficients of ${\bf U}_k(z,s)$ coincide with those of ${\bf
F}_k(z)$ and are independent of $s$.\qed

We now consider the random variable $X_{k,\tau}^n$ having the
distribution
\begin{equation}\label{E:sigma}
\mathbb{P}(X_{k,\tau}^n=t)= {\sf T}_{k,\tau}(n,t)/{\sf
T}_{k,\tau}(n) \,
\end{equation}
where $t=0,1,\ldots \lfloor {n}/{2}\rfloor$. We shall show in the 
following theorem that the distribution of $X_{k,\tau}^n$ is determined by
the shift of the singularity parametrized by $s$. In other words,
Theorem~\ref{T:z} induces the below central limit theorem. The
proof-idea of Theorem~\ref{T:Gauss} is due to Bender, where it
appeared in slightly different context \cite{T:Gauss}. It follows
from analyzing the characteristic function, using the
{L\'{e}vy-Cram\'{e}r} Theorem of Section~\ref{S:basic}.
\begin{theorem}\label{T:Gauss}
For any $2\le k\le 9$ and $3\le \tau\le 7$ there exist a pair
$(\mu_{k,\tau},\sigma_{k,\tau})$ such that the normalized random
variable
\begin{equation}
Y_{k,\tau}^n=\frac{X_{k,\tau}^n- \mu_{k,\tau} \, n}{\sqrt{{n\,
\sigma_{k,\tau}}^2 }}
\end{equation}
has asymptotically normal distribution with parameter $(0,1)$. That
is we have
\begin{equation}\label{E:converge}
\lim_{n\to\infty}\mathbb{P}\left(\frac{X_{k,\tau}^n- \mu_{k,\tau}
n}{\sqrt{n\, \sigma_{k,\tau}^2}} < x \right)  =
\frac{1}{\sqrt{2\pi}}\int_{-\infty}^{x}\,e^{-\frac{1}{2}c^2} dc \ ,
\end{equation}
where $\mu_{k,\tau}$ and $\sigma_{k,\tau}^2$ are given by
\begin{equation}\label{E:was}
\mu_{k,\tau}= -\frac{\gamma_{k,\tau}'(0)}{\gamma_{k,\tau}(0)},
\qquad \qquad \sigma_{k,\tau}^2=
\left(\frac{\gamma_{k,\tau}'(0)}{\gamma_{k,\tau}(0)}
\right)^2-\frac{\gamma_{k,\tau}''(0)}{\gamma_{k,\tau}(0)}.
\end{equation}
\end{theorem}


\pff Suppose we are given the random variable (r.v.) $\xi_n$ with
mean $\mu_n$ and variance $\sigma_n^2$. We consider the rescaled
r.v.~$\eta_{n}=(\xi_{n} - \mu_{n})\sigma^{-1}_{n}$ and the
characteristic function of $\eta_{n}$:
\begin{equation}
f_{\eta_n}(c)=\mathbb{E}[e^{ic\eta_{n}}]=
 \mathbb{E}[e^{ic\frac{\xi_{n}}{\sigma_{n}}}]
                                     e^{-i\frac{\mu_{n}}{\sigma_{n}}c} \ .
\end{equation}
Writing $X_n$ instead of $X_{k,\tau}^n$ we derive for $\xi_n=X_n$,
substituting for the term $\mathbb{E}[e^{ic\eta_{n}}]$
\begin{equation}
f_{X_n}(c)=\left( \sum_{t=0}^{n/2}\frac{{\sf T}_{k,\tau}(n,t)} {{\sf
T}_{k,\tau}(n)} e^{ic\frac{t}{\sigma_{n}}}\right)\,
e^{-i\frac{\mu_{n}}{\sigma_{n}}c}
 \ .
\end{equation}
In view of $\varphi_{n,k,\tau}(s)= \sum_{t\le n/2}{\sf
T}_{k,\tau}(n,t)e^{ts}$ we interpret
$$ \varphi_{n,k,\tau}(0)=\sum_{t\le n/2}{\sf T}_{k,\tau}(n,t)
\quad\text{\rm  and}\quad \varphi_{n,k,\tau}(({ic})/({\sigma_n}))=
\sum_{t\le n/2}{\sf T}_{k,\tau}(n,t)e^{t({ic})/({\sigma_n})} \ .
$$
Writing $\varphi_n$ instead of $\varphi_{n,k,\tau}$, we accordingly
obtain
\begin{equation}\label{E:well}
f_{X_n}(c)=\frac{1}{\varphi_{n}(0)}\,\varphi_{n}\left(\frac{ic}{\sigma_{n}}
\right)\, e^{-i\frac{\mu_{n}}{\sigma_{n}}c} \ .
\end{equation}
Now we have to provide the interpretation of $\varphi_n(0)$ and
$\varphi_n(({ic})/({\sigma_n}))$. This is facilitated via
Theorem~\ref{T:z}:
\begin{equation}\label{E:4}
[z^n]\,{\bf U}_k(z,s)= K(s) \,\theta_k(n)\,
\left(\gamma_{k,\tau}(s)^{-1}\right)^n
\left(1-O({1}/{n})\right)\quad \text{\rm for some
$K(s)\in\mathbb{C}$},
\end{equation}
uniformly in $s$ and where $\theta_k(n)$ is some subexponential
factor, independent of $s$ (we showed that the singular expansion
remains invariant when substituting $\psi_{\tau}(z,s)$ for $z$).
Therefore
\begin{equation}\label{E:uu}
f_{X_n}(c)\sim \frac{K(\frac{ic}{\sigma_n})}{K(0)}\,
\left[\frac{\gamma_{k,\tau}(\frac{ic}{\sigma_n})}
{\gamma_{k,\tau}(0)}\right]^{-n} e^{-i\frac{\mu_{n}}{\sigma_{n}}c},
\end{equation}
uniformly in $c$, where $c$ is contained in an arbitrary bounded
interval. Taking the logarithm we obtain
\begin{equation}
\ln f_{X_n}(c)\sim\ln\frac{K(\frac{ic}{\sigma_n})}{K(0)}- n \,
\ln\frac{\gamma_{k,\tau}(\frac{ic}{\sigma_n})} {\gamma_{k,\tau}(0)}
-i \frac{\mu_n}{\sigma_n}c \ .
\end{equation}
Expanding $g(s)=\ln({\gamma_{k,\tau}(s)})/( {\gamma_{k,\tau}(0)})$
in its Taylor series at $s=0$, (note that $g(0)=0$ holds) yields
\begin{equation}\label{E:kkk}
\ln\frac{\gamma_{k,\tau}(\frac{ic}{\sigma_{n}})}
{\gamma_{k,\tau}(0)}= \frac{\gamma_{k,\tau}'(0)}{\gamma_{k,\tau}(0)}
\frac{ic}{\sigma_{n}}-
\left[\frac{\gamma_{k,\tau}''(0)}{\gamma_{k,\tau}(0)}-
\left(\frac{\gamma_{k,\tau}'(0)}
{\gamma_{k,\tau}(0)}\right)^2\right] \frac{c^2}{2\sigma^2_{n}}
+{O}(\left(\frac{ic}{\sigma_n}\right)^{3})
\end{equation}
and $\ln f_{X_n}(c)$ becomes asymptotically
\begin{equation}\label{E:nun}
\begin{split}
& \ln \frac{K(\frac{ic}{\sigma_n})}{K(0)}  - n \,
\left\{\frac{\gamma_{k,\tau}'(0)}
{\gamma_{k,\tau}(0)}\frac{ic}{\sigma_{n}}-\frac{1}{2}
\left[\frac{\gamma_{k,\tau}''(0)}{\gamma_{k,\tau}(0)}-
\left(\frac{\gamma_{k,\tau}'(0)}{\gamma_{k,\tau}(0)}
\right)^2\right] \frac{c^2}{\sigma^2_{n}}
+{O}(\left(\frac{ic}{\sigma_n}\right)^{3})\right\}\\
&  -
\frac{i\mu_{n}c}{\sigma_{n}}.
\end{split}
\end{equation}
${\bf U}_k(z,s)$ is analytic in $s$ where $s$ is contained in a disc
of radius $\epsilon$ around $0$ and therefore in particular
continuous in $s$ for $\vert s\vert <\epsilon$. In view of
eq.~(\ref{E:nun}) we set
$$
\mu=-\frac{\gamma_{k,\tau}'(0)}{\gamma_{k,\tau}(0)}, \quad \quad
\sigma^2=
\left(\frac{\gamma_{k,\tau}'(0)}{\gamma_{k,\tau}(0)}
\right)^2-\frac{\gamma_{k,\tau}''(0)} {\gamma_{k,\tau}(0)}
$$
Setting $\mu_n=n \mu$ and $\sigma_n^2=n\sigma^2$ we can conclude
from eq.~(\ref{E:4}) for fixed $c\in ]-\infty, \infty[$
\begin{equation}
\lim_{n\to \infty}\left(\ln K(({ic})/({\sigma_n}))-\ln K(0)\right)=0
\end{equation}
and eq.~(\ref{E:nun}) becomes
\begin{equation}
\ln f_{X_n}(c)\sim -{\; c^2}/{2}
+{O}(\left(({ic})/{\sigma_n}\right)^{3})
\end{equation}
with uniform error term for $c$ from any bounded interval. This is
equivalent to $\lim_{n\rightarrow\infty}f_{X_n}(c)=
\exp(-{c^2}/{2})$  uniformly in $c$. The L\'{e}vy-Cram\'{e}r Theorem
of Section~\ref{S:basic} implies now eq.~(\ref{E:converge}) and the
proof of Theorem~\ref{T:Gauss} is complete.\qed

{\bf Acknowledgments.}
We are grateful to Gang Ma for his comments.
This work was supported by the 973 Project, the PCSIRT Project of
the Ministry of Education, the Ministry of Science and Technology,
and the National Science Foundation of China.


\bibliographystyle{plain}


\newpage

\begin{table}
\begin{center}
\begin{tabular}{|ccccccc|cc|}
\hline
&\multicolumn{2}{c}{$k=2$} &\multicolumn{2}{c}{$k=3$}        & \multicolumn{2}{c|}{$k=4$} \\
\hline
& $\mu_{k,\tau}$ & $\sigma_{k,\tau}^2$ & $\mu_{k,\tau}$ & $\sigma_{k,\tau}^2$ & $\mu_{k,\tau}$ & $\sigma_{k,\tau}^2$ \\
\hline
$\tau=3$ & 0.090323 & 0.0189975 & 0.115473 & 0.0086760 & 0.123509 & 0.0076977 \\
$\tau=4$ & 0.071677 & 0.0131316 & 0.086554 & 0.0055685 & 0.091737 & 0.0049917 \\
$\tau=5$ & 0.059591 & 0.0098165 & 0.069467 & 0.0039688 & 0.073166 & 0.0035769 \\
$\tau=6$ & 0.051092 & 0.0077233 & 0.058149 & 0.0026885 & 0.060964 & 0.0027313 \\
$\tau=7$ & 0.044774 & 0.0062991 & 0.050083 & 0.0017584 & 0.052319 & 0.0021788\\
\hline
&\multicolumn{2}{c}{$k=5$}  &\multicolumn{2}{c}{$k=6$}                        & \multicolumn{2}{c|}{$k=7$} \\
\hline
 & $\mu_{k,\tau}$ & $\sigma_{k,\tau}^2$ & $\mu_{k,\tau}$ & $\sigma_{k,\tau}^2$ & $\mu_{k,\tau}$ & $\sigma_{k,\tau}^2$ \\
\hline
$\tau=3$ & 0.128157 & 0.0070210 & 0.131353 & 0.0065187 & 0.133748 & 0.0061254 \\
$\tau=4$ & 0.094768 & 0.0020037 & 0.119551 & 0.0080515 & 0.098461 & 0.0040797 \\
$\tau=5$ & 0.075345 & 0.0033114 & 0.076864 & 0.0031162 & 0.078016 & 0.0029639 \\
$\tau=6$ & 0.062629 & 0.0025364 & 0.063794 & 0.0023936 & 0.064680 & 0.0022823 \\
$\tau=7$ & 0.053648 & 0.0020277 & 0.054580 & 0.0019171 & 0.055291 & 0.0018310\\
\hline
\end{tabular}
\centerline{}
\smallskip
\caption{\small The mean $\mu_{k,\tau}$ and variance $\sigma_{k,\tau}^2$
of the normal limit distributions of the random variable $X_{k,\tau}^n$, for 
different $k$ and $\tau$.} \label{T:1}
\end{center}
\end{table}

\newpage

\newpage

\newpage

\newpage

\end{document}